\documentclass[a4paper, 11pt]{article}
\usepackage{fullpage}
\usepackage{amsmath}
\usepackage{amsthm}
\usepackage{amsfonts}
\usepackage{amssymb}
\usepackage{tikz-cd}
\usepackage[utf8x]{inputenc}
\usepackage{enumerate}
\usepackage[english]{babel}
\usepackage{esint}
\usepackage{hyperref}
\usepackage{bm}
\usepackage{bbm}
\usepackage[mathscr]{eucal}

\usepackage{soul}
\allowdisplaybreaks[1]

\newtheorem{Theorem}{Theorem}[section]
\newtheorem{Proposition}[Theorem]{Proposition}
\newtheorem{Lemma}[Theorem]{Lemma}
\newtheorem{Corollary}[Theorem]{Corollary}

\theoremstyle{definition}
\newtheorem{Definition}[Theorem]{Definition}

\newtheorem{Remark}[Theorem]{Remark}

\newcommand{\bTheorem}[1]{
\begin{Theorem} \label{T#1} }
\newcommand{\eT}{\end{Theorem}}

\newcommand{\bProposition}[1]{
\begin{Proposition} \label{P#1}}
\newcommand{\eP}{\end{Proposition}}

\newcommand{\bLemma}[1]{
\begin{Lemma} \label{L#1} }
\newcommand{\eL}{\end{Lemma}}

\newcommand{\bCorollary}[1]{
\begin{Corollary} \label{C#1} }
\newcommand{\eC}{\end{Corollary}}

\newcommand{\bRemark}[1]{
\begin{Remark} \label{R#1} }
\newcommand{\eR}{\end{Remark}}

\newcommand{\bDefinition}[1]{
\begin{Definition} \label{D#1} }
\newcommand{\eD}{\end{Definition}}

\newcommand{\Del}{\Delta_x}

\newcommand{\Ds}{\mathbb{D}_x}

\newcommand{\bfphi}{\boldsymbol{\varphi}}

\newcommand{\bFormula}[1]{
\begin{equation} \label{#1}}
\newcommand{\eF}{\end{equation}}

\newcommand{\Ov}[1]{\overline{#1}}

\newcommand{\DC}{C^\infty_c}

\newcommand{\vr}{\varrho}
\newcommand{\vre}{\vr_\ep}

\newcommand{\vue}{\vu_\ep}

\newcommand{\tvu}{{\tilde \vu}}

\newcommand{\vu}{\vc{u}}
\newcommand{\vm}{\vc{m}}

\newcommand{\vn}{\vc{n}}

\newcommand{\vc}[1]{{\bf #1}}

\newcommand{\Div}{{\rm div}_x}
\newcommand{\Grad}{\nabla_x}

\newcommand{\dx}{\,{\rm d} {x}}

\newcommand{\dt}{\,{\rm d} t }

\newcommand{\dxdt}{\dx  \dt}
\newcommand{\intO}[1]{\int_{\Omega} #1 \ \dx}

\newcommand{\intTO}[1]{\int_0^T \int_{\Omega} #1 \ \dxdt}

\newcommand{\vv}{\vc{v}}

\newcommand{\D}{{\rm d}}

\newcommand{\ep}{\varepsilon}

\renewcommand{\S}{\mathbb{S}}

\newcommand{\am}{{\alpha_m}}

\def\softd{{\leavevmode\setbox1=\hbox{d}%
          \hbox to 1.05\wd1{d\kern-0.4ex{\char039}\hss}}}
\definecolor{Cgrey}{rgb}{0.85,0.85,0.85}
\definecolor{Cblue}{rgb}{0.50,0.85,0.85}
\definecolor{Cred}{rgb}{1,0,0}
\definecolor{fancy}{rgb}{0.10,0.85,0.10}

\newcommand\Cbox[2]{%
    \newbox\contentbox%
    \newbox\bkgdbox%
    \setbox\contentbox\hbox to \hsize{%
        \vtop{
            \kern\columnsep
            \hbox to \hsize{%
                \kern\columnsep%
                \advance\hsize by -2\columnsep%
                \setlength{\textwidth}{\hsize}%
                \vbox{
                    \parskip=\baselineskip
                    \parindent=0bp
                    #2
                }%
                \kern\columnsep%
            }%
            \kern\columnsep%
        }%
    }%
    \setbox\bkgdbox\vbox{
        \color{#1}
        \hrule width  \wd\contentbox %
               height \ht\contentbox %
               depth  \dp\contentbox
        \color{black}
    }%
    \wd\bkgdbox=0bp%
    \vbox{\hbox to \hsize{\box\bkgdbox\box\contentbox}}%
    \vskip\baselineskip%
}

\numberwithin{equation}{section}

\date{}
\begin{document}



\title{Generalized dissipative solutions to free boundary \\compressible viscous models}

\author{Anna Abbatiello \and  Donatella Donatelli}


\maketitle

\centerline{Department of Mathematics and Physics, University of Campania ``L.~Vanvitelli"}
\centerline{Viale A.~Lincoln 5, 81100 Caserta, Italy}

\centerline{anna.abbatiello@unicampania.it}

\centerline{Department of Information Engineering, Computer Science and Mathematics}
\centerline{University of L'Aquila, 67100 L'Aquila, Italy}

{\centerline{donatella.donatelli@univaq.it}}

\begin{abstract} 
We study free boundary compressible viscous models that may include nonlinear viscosities. These are compressible/incompressible Navier-Stokes type systems for a non-Newtonian stress tensor. They describe the motion of a possibly non-Newtonian fluid in free flow and in congested regions. In the congested regions it appears the pressure that is the Lagrange multiplier associated with the incompressibility constraint, while in free flows it is a pressureless gas system. 
We establish the existence of generalized dissipative solutions in the case of in/out-flow boundary conditions and we also prove that if these solutions are  smooth then they are  classical solutions.   
\end{abstract}


\section{Introduction}
\label{P}
In this paper we are concerned with the following compressible/incompressible Navier-Stokes type of system
\begin{subequations} \label{problem}\begin{align}
\partial_t \vr + \Div (\vr \vu) &= 0  \ \mbox{ in } (0, T)\times\Omega, \label{c-e}\\
\partial_t (\vr \vu) + \Div (\vr \vu \otimes \vu) -  \Div \mathbb{S}+ \Grad \pi &= 0  \ \mbox{ in } (0, T)\times\Omega,\label{mom-e}\\
0\leq \vr &\leq 1 \ \mbox{ in } (0, T)\times\Omega,\label{bounded}\\
\Div \vu &= 0 \mbox{ on } \{ \vr=1\}, \label{r=1}\\
\pi=0 \mbox{ in } \{\vr<1\}, \ \pi&\geq 0 \mbox{  on } \{\vr=1\},\label{pi=0}
\end{align}
\end{subequations}
with $T>0$ and $\Omega \subset \mathbb{R}^d$  a bounded domain of dimension $d=3$.
The unknowns of the system are $\vr=\vr(t,x)$  the density,  $\vu=\vu(t,x)$  the velocity and $\pi$ the pressure term,  $\mathbb{S}$ is the viscous stress tensor. The system  \eqref{problem} describes a free boundary problem for a  fluid, defined by means of a pressure threshold above which the fluid is taken to be incompressible, and below which the fluid is compressible and was introduced by Lions and  Masmoudi in \cite{LiMa}.  The  motivation to study this free boundary problems  comes from a modelling point of view. This set of equations can be used  to  study  fluids with large bubbles filled with gas, those models involve a threshold on the pressure beyond which one has the incompressible Navier-Stokes equations for the fluid and below which one has a compressible model for the gas. Another motivation comes from collective dynamics for the motion of large crowds, where $\vr=1$ represents the congested density. This means that when the density of the crowd reaches the value $\vr=1$, the crowd behaves like an incompressible fluid. When the density is less than $1$ the particles move freely and  the crowd behaves like a compressible fluid, see for examples the paper \cite{DMZ18}. Other applications  of free boundary models can be find in 
\cite{DT2017}, \cite{DT2018}, \cite{DTT2020} in the case of polymeric fluids or in  \cite{PZ15}, \cite{PoWrZa}  for two-phase flows.

The purpose of this paper is to study a twofold generalization of the model \eqref{problem}. 

First we will consider the fluid in the case of non-Newtonian stress tensor. Indeed, the viscous stress tensor $\mathbb{S}$ is related to the 
symmetric velocity gradient 
$\Ds \vu = \frac{1}{2}\left(\Grad \vu + (\Grad \vu)^t\right)$
through a general  implicit rheological law
\begin{equation} \label{P6}
\Ds \vu : \mathbb{S} = 
F(\Ds \vu) + F^*(\mathbb{S}) 
\end{equation}
for a given suitable convex potential $F$ and its conjugate $F^*$ defined in the Section \ref{Viscous stress}.  To the best of our knowledge the  non-Newtonian context for this kind of free boundary problems are almost unexplored. 

The second generalization concerns the fact that we will study the problem in the framework of open fluids systems that interact with the surroundings.
From a modelling point of view an open system is described by a proper choice of boundary conditions. An example of boundary conditions in the case of a closed system is given by the no-slip boundary conditions. In the case of open system the no-slip boundary condition is replaced by a  general Dirichlet boundary condition
$$\vu|_{\partial \Omega} = \vu_B,$$
where $ \vu_B= \vu_B(x)$ is a given field that we assume depends on the space variable and satisfies certain properties, see Section \ref{boundary data}. For a complete overview of compressible open fluid systems see \cite{FNopen} and references therein.

The main goal of the present paper is to develop a mathematical theory for this free boundary model
{in the case of complicated material responses that may be expressed in a compact way through an \emph{implicit relation}.  We formulate the material response in terms of the potential of the viscous stress tensor and we only assume a growth condition. We include for example the general case of the power-law viscous fluids for which the existence of standard weak solution is not known in the compressible case.  In a such general framework the corresponding energy estimate provides only bounds of the total variation of the velocity. For this reason we rewrite the energy inequality using the so-called Fenchel-Young inequality which allows the limit passage thanks to the convexity of the potential and of its convex conjugate. }

We can prove the existence of the so called \emph{dissipative generalized solutions}, those are solutions that satisfy the standard weak formulation modulo a defect measure, for related results in this framework see  \cite{AbFeNo}.  The main  features of the dissipative weak solutions is that the energy inequality (balance) is an integral part of the definition of the dissipative solutions and the family of unknown is augmented  by the  quantity 
$$\mathfrak{R} \in L^\infty(0,T; \mathcal{M}^+(\Ov{\Omega}; R^{d \times d}_{\rm sym}))$$
called the Reynolds viscous stress that accounts e.g. for possible concentrations in the convective term. We end this section with a plan of the paper. In Section \ref{formulation} we set up our problem, give the definition of dissipative weak solutions and state our main results. In Section \ref{E} we describe the multilevel approximation scheme needed to prove our existence results. Section \ref{existence} is devoted to the proof of the existence of dissipative solutions and finally in Section \ref{compatibility} we show that dissipative solutions enjoying additional regularity properties are classical solutions.

\section{Formulation of the problem and main results}
\label{formulation}

\subsection{Viscous stress}
 \label{Viscous stress}
 We are interested in non-Newtonian fluid-like responses. One prototype example is a potential of type 
 $$ F(\mathbb{D})=\frac{\mu_0}{q}(\mu_1 + |\mathbb{D}_0|^2)^{\frac{q}{2}} + \eta_0 (\eta_1+ \left| {\rm tr}[\mathbb{D}] \right|^2)^{\frac{q}{2}} \mbox{ with } \mu_0>0, \ \mu_1 \geq 0, \ \eta_0\geq0, \ \eta_1\geq0,\ q>1, $$
where $$\mathbb{D}_0=\mathbb{D} - \frac{1}{d} {\rm tr}[\mathbb{D}] \mathbb{I}.$$
The relation between the symmetric velocity gradient $\Ds \vu$ and the viscous stress $\mathbb{S}$ is determined by the choice 
of the potential $F$. We assume that 
\begin{equation} \label{S3}
F: \mathbb{R}^{d \times d}_{\rm sym} \to [0, \infty) \ \mbox{is a (proper) convex function},\ F(0) = 0,
\end{equation}
and that there exist $\mu > 0$ and $q > 1$ such that
\begin{equation} \label{S5a}
F(\mathbb{D}) \geq \mu \left| \mathbb{D} - \frac{1}{d} {\rm tr}[\mathbb{D}] \mathbb{I} \right|^q
\ \mbox{for all}\ |\mathbb{D}| > 1.
\end{equation}

\subsection{Initial Data}
The initial conditions are as follows:
\begin{equation} \label{P7}
\vr(0, \cdot) = \vr_0,  \ \vr \vu(0, \cdot) = \vm_0,
\end{equation}
where
\begin{eqnarray}
& & 0 \le \vr_{0} \le 1, \quad  \vr_{0} \in L^1(\Omega), \quad \vr_{0} \not \equiv 0,\,\,\, \int_{\Omega}\!\!\!\!\!\!\!-\vr_{0}=M<1,\nonumber\\
& & \vm_{0} \in L^2 (\Omega), \quad \vm_{0} =0 \,\, \mbox{a.e.  on }\,\, \{\vr_{0} = 0\},\nonumber\\
&& \vr_{0} |\vu_{0}|^{2}  \in L^2(\Omega),\quad \vu^0 = \frac{\vm_{0}}{\vr_{0}}\,\, \mbox{on} \,\, \{\vr_{0} >0\}\nonumber. 
\end{eqnarray}

\subsection{Boundary  Data}
\label{boundary data}
 We define the following subsets of $\Omega$ 
 $$\Gamma_{in} = \left\{ x \in \partial \Omega \ \Big| \ \vu_B \cdot \vc{n} < 0 \right\}, \quad 
 \Gamma_{out} = \left\{ x \in \partial \Omega \ \Big| \ \vu_B \cdot \vc{n} \geq 0 \right\},$$ 
 where $\vc{n}$ is the outer normal vector to $\partial \Omega$.
 We prescribe general inflow--outflow boundary conditions.
 \begin{align} \label{P5}
\vu|_{\partial \Omega} &= \vu_B,\\
\vu_B&\in C^2(\partial\Omega;\mathbb{R}^d), \quad  \int_{\partial\Omega} \vu_B\cdot\vn \geq 0, \label{bdu}\\
 \vr|_{\Gamma_{in}} &= \vr_B, \quad  \vr_B \in C^1(\partial \Omega), \quad  0<\vr_B\leq 1  \mbox{ on } \Gamma_{\rm in}, \quad \inf_{\Gamma_{\rm in}} \vr_B>0. \label{hyp-rb}
\end{align}

\begin{Remark}
In both our results we need to assume the velocity flux has a sign
$$\int_{\partial\Omega} \vu_B\cdot\vn \geq 0.$$
This assumption is needed when we pass into the limit in Galerkin approximations. In this last step we need to extend $\vu_{B}$ to the whole domain $\Omega$ in such a way that this extension is smooth and $\Div \vu_B \geq 0 $. In the case of a negative flux we are not able to get uniform bounds by the energy estimate.
\end{Remark}
\subsection{Definition of generalized dissipative solutions} 
We introduce the notion of generalized dissipative solution to system \eqref{problem}.

\begin{Definition}\label{def}
A quadruplet $(\vr, \vu, \S, \pi)$ is a generalized dissipative solution to system \eqref{problem} if the following conditions are satisfied.
\begin{enumerate}
\item[(i)] The  density $\vr$ satisfies
 $$0\leq \vr\leq 1$$ 
 and belongs to the class
$$\vr \in C_{\rm weak}([0,T]; L^\gamma(\Omega)) \cap L^\gamma(0,T; L^\gamma(\partial \Omega; |\vu_B \cdot \vc{n}|
{\rm d}S_x)),\
\mbox{ for any }\gamma > 1. $$
The velocity $\vu$ belongs to the class
$$\vu \in L^q (0,T; W^{1,q}(\Omega; R^d)),\ (\vu - \vu_B) \in L^q(0,T; W^{1,q}_0(\Omega; R^d)),\ q > 1,$$
also it holds that $\vm \equiv \vr \vu$. There exist
\[
\mathbb{S} \in L^1((0,T) \times \Omega; R^{d \times d}_{\rm sym}),\ \mathfrak{R} \in L^\infty(0,T; \mathcal{M}^+(\Ov{\Omega}; R^{d \times d}_{\rm sym})),
\]
where the symbol
$\mathcal{M}^+(\Ov{\Omega}; R^{d \times d}_{\rm sym})$ denote the set of all positively semi--definite tensor valued measures on $\Ov{\Omega}$. 
The momentum satisfies
\[
\vr \vu \in C_{\rm weak}([0,T]; L^{\frac{2 \gamma}{\gamma + 1}}(\Omega; R^d)).
\]

The pressure $$ \pi \in \mathcal{M}^+(I\times K) \mbox{ for any compact }I\times K\subset (0,T)\times \Omega$$
satisfies 
\begin{equation*}
\begin{split}
&(\vr-1)\pi=0 \mbox{ a.e. in } (0,T)\times\Omega,\\
&\Div \vu=0 \mbox{ a.e. on } \{\vr=1\}.
\end{split}\end{equation*}
\item [(ii)] The continuity equation
\begin{equation*}
\begin{split}
\left[ \intO{ \vr \varphi } \right]_{t = 0}^{t = T}&= 
\int_0^T \intO{ \Big[ \vr \partial_t \varphi + \vr \vu \cdot \Grad \varphi  \Big] }
\dt \\ 
&- \int_0^T \int_{\Gamma_{\rm out}} \varphi \vr \vu_B \cdot \vc{n} \ {\rm d} S_x \dt - 
\int_0^T \int_{\Gamma_{\rm in}} \varphi \vr_B \vu_B \cdot \vc{n}  \ {\rm d}S_x \dt,
\end{split}
\end{equation*}
holds for any $\varphi \in C^1([0,T] \times \Ov{\Omega})$. Also $\vr(0,\cdot)=\vr_0$.
\item[(iii)] 
The momentum  equation 
\begin{equation*} 
\begin{split}
\left[ \intO{ \vr \vu \cdot \bfphi } \right]_{t=0}^{t = T} &= 
\int_0^T \intO{ \Big[ \vr \vu \cdot \partial_t \bfphi + \vr \vu \otimes \vu : \Grad \bfphi 
 - \mathbb{S}  : \Grad \bfphi \Big] } \\
 &+\langle \pi, \Div\bfphi\rangle_{(\mathcal{M}((0,T)\times\Omega), C((0,T)\times\Omega))} + \int_0^T\int_\Omega\Grad\bfphi:\D\mathfrak{R}(t)\dt
\end{split}
\end{equation*}
holds for any $\bfphi \in C_c^1((0,T)\times\Omega; \mathbb{R}^d)$.
\item[(iv)] The energy inequality holds
\begin{equation*} 
\begin{split}
\left[ \intO{\frac{1}{2} \vr |\vu - \vu_B|^2 } \right]_{t = 0}^{ t = T} &+ 
\int_0^T \intO{\left[ F(\Ds \vu) + F^*(\mathbb{S}) \right] } \dt 
+\langle \pi, \Div\vu_B\rangle_{(\mathcal{M}((0,T)\times\Omega), C((0,T)\times\Omega))}\\
   + \frac{1}{\Ov{d}} \int_{\Ov{\Omega}}  \D \ {\rm tr}[\mathfrak{R}] (T)
&\leq -\int_0^T \int_{\Ov{\Omega}} \Grad \vu_B : \D \ \mathfrak{R}(t) \dt - 
\int_0^T \intO{ \left[ \vr \vu \otimes \vu \right] :  \Grad \vu_B } \dt \\ &+ \int_0^T \intO{ {\vr} \vu  \cdot \vu_B \cdot \Grad \vu_B  } 
\dt
+ \int_0^T \intO{ \mathbb{S} : \Grad \vu_B } \dt
\end{split}
\end{equation*}
with $\Ov{d}$ a constant depending on the dimension $d$.
\end{enumerate}
\end{Definition}

\begin{Remark}
We point out that the condition $0\leq \vr\leq 1$ has to be understood in the sense of almost everywhere defined functions and since $\pi$ is not a function  the condition \eqref{pi=0} makes sense in the equivalent formulation $(\vr-1)\pi=0$ that we have in (i) of the above definition.
\end{Remark}

\begin{Remark}
It is not difficult to show that if $\mathfrak{R}=0$, then $\vr$ and $\vu$ are standard weak solutions of \eqref{problem} fulfilling a variant of the energy inequality.
\end{Remark}

\subsection{Main results}
Now we can state the main result of this paper. The first one concerns the existence of solutions to problem \eqref{problem} and reads as follows.
\begin{Theorem}
Let the assumptions \eqref{P7}--\eqref{hyp-rb}  be satisfied and let $T>0$ then problem \eqref{problem} admits at least one generalized dissipative solution in the sense of Definition \ref{def}.
\label{T1}
\end{Theorem}

The second results is related to the compatibility between generalized dissipative and strong solutions. We show that if a dissipative solution satisfies additional regularity properties then it is a classical solutions (i.e. $\mathfrak{R}=0$ and $\S\in \partial F(\mathbb{D})$).

\begin{Theorem}
Let the assumptions \eqref{P7}--\eqref{hyp-rb} be satisfied, let $T>0$ and let $(\vr, \vu, \S, \pi)$ be a generalized dissipative solution to system \eqref{problem} in the sense of Definition \ref{def} belonging to the class
\begin{equation*}
\vu\in C^1([0,T]\times \Ov{\Omega}; \mathbb{R}^d), \ \vr\in C^1([0,T]\times \Ov{\Omega}), \inf_{(0,T)\times\Omega}\vr >0, \ \pi\in C([0,T]\times\Ov{\Omega}),
\end{equation*}
then the equations are satisfied in the classical sense and $\mathfrak{R}=0$ and $\S\in \partial F(\mathbb{D})$.
\label{T2}
\end{Theorem}

\section{The multilevel approximating problem}
\label{E}
The proof of the Theorem \ref{T1} follows from a multilevel approximation which combines some features of the compressible Navier Stokes equations accommodated in such a way to deal with the non-Newtonian viscosity term and the  free boundary problem. 
First, we introduce a sequence of finite--dimensional spaces 
$X_n \subset L^2(\Omega; \mathbb{R}^d)$,
\[
X_n = {\rm span} \left\{ \vc{w}_i\ \Big|\ \vc{w}_i \in \DC(\Omega; R^3),\ i = 1,\dots, n \right\}
\]   
such that the sequence $\{\vc{w}_i\}_i$ is made of orthonormal functions with respect to the standard scalar product in $L^2$. Then we regularize the convex potential $F$ and we perform a parabolic approximation of the continuity equation combined with a Galerkin approximation of the momentum equation. Finally, in order to deal with the free boundary problem we introduce a sequence of approximating isoentropic pressures with the adiabatic exponent that goes to infinity. 

\subsection{Regularization of the potential}
We regularize the convex potential $F$ to make it continuously differentiable: 
\[
F_\delta(\mathbb{D}) = \int_{R^{d \times d}_{\rm sym}} \xi_\delta ( | \mathbb{D} - \mathbb{Z} | ) F(\mathbb{Z}) \ {\rm d} \mathbb{Z}
- \int_{R^{d \times d}_{\rm sym}} \xi_\delta ( | \mathbb{Z} | ) F(\mathbb{Z}) \ {\rm d} \mathbb{Z}
\]
where $ \{\xi_\delta\}$ with $\delta>0$ is a family 
of regularizing kernels. The regularized potentials $F_\delta$ are convex, non--negative, infinitely differentiable, $F_\delta (0) = 0$, and satisfy \eqref{S5a}, specifically
\begin{equation} \label{F1a}
F_\delta (\mathbb{D}) \geq \mu \left| \mathbb{D} - \frac{1}{d} {\rm tr}[\mathbb{D}] \mathbb{I} \right|^q
\ \mbox{for all}\ |\mathbb{D}| > 1,
\end{equation}
with $q > 1$, $\mu > 0$ independent of $\delta \searrow 0$.
\subsection{Parabolic approximation of the equation of continuity}
We introduce a parabolic approximation of the equation of continuity, namely  for any $\ep>0$ we set
\begin{subequations}\label{Den}\begin{align} \label{E1}
&\partial_t \vr + \Div (\vr \vu ) = \ep \Del \vr \ &&\mbox{in}\ (0,T) \times \Omega,\\
 \label{E2}
&\ep \Grad \vr \cdot \vc{n} + (\vr_B - \vr) [\vu_B \cdot \vc{n}]^- = 0 \ &&\mbox{on}\ (0,T) \times \partial \Omega,\\
 \label{E3}
&\vr(0, \cdot) = \vr_0 \ && \mbox{in } \Omega,
\end{align}\end{subequations}
where $\vu = \vv + \vu_B$, with $\vv \in C([0,T]; X_n)$, $\vu|_{\partial\Omega}=\vu_B$. Given $\vu$, $\vr_B$, $\vu_B$, this is a linear problem for the unknown $\vr$. We avoid to assume regularity of  $\Omega$ and suppose that it is Lipschitz, therefore we introduce the weak formulation of problem \eqref{Den}
\begin{equation} \label{E4}
\begin{split}
&\left[ \intO{ \vr \varphi } \right]_{t = 0}^{t = \tau}= 
\int_0^\tau \intO{ \left[ \vr \partial_t \varphi + \vr \vu \cdot \Grad \varphi - \ep \Grad \vr \cdot \Grad \varphi \right] }
\dt \\ 
&- \int_0^\tau \int_{\partial \Omega} \varphi \vr \vu_B \cdot \vc{n} \ {\rm d} S_x \dt + 
\int_0^\tau \int_{\partial \Omega} \varphi (\vr - \vr_B) [\vu_B \cdot \vc{n}]^{-}  \ {\rm d}S_x \dt,\ 
\vr(0, \cdot) = \vr_0, 
\end{split}
\end{equation}
for any test function
$\varphi \in L^2(0,T; W^{1,2}(\Omega))$, $\partial_t \varphi \in L^1(0,T; L^2(\Omega))$. 

 By virtue of the existence result by Crippa-Donadello-Spinolo \cite[Lemma 3.2]{CrDoSp}, we claim that given $\vu = \vv + \vu_B$,  the initial--boundary value problem \eqref{E1}--\eqref{E3} 
admits a weak solution $\vr$, unique in the class 
\[
\vr \in L^2(0,T; W^{1,2}(\Omega)) \cap C([0,T]; L^2(\Omega)).
\]
The norm in the aforementioned spaces is bounded only in terms of the data $\vr_B$, $\vr_0$, $\vu_B$ and 
\[
\sup_{t \in (0,T)}
\| \vv(t, \cdot) \|_{X_n}.
\]
Moreover it holds the maximum principle (see \cite[Lemma 3.4]{CrDoSp}) 
\begin{equation}\label{max}
\| \vr \|_{L^\infty((0,\tau) \times \Omega)} \leq 
\max \left\{ \| \vr_0 \|_{L^\infty(\Omega)}; \| \vr_B \|_{L^\infty((0,T) \times \Gamma_{in})}; 
\| \vu_B \|_{L^\infty((0,T) \times \Omega)} \right\} \exp \left( \tau \| \Div \vu \|_{L^\infty((0,\tau) \times \Omega)} \right).
\end{equation}
We can perform renormalization of equation \eqref{E1}, i.e. let $B \in C^2(R)$ then the renormalized equation
\begin{equation}\label{renorm}
\begin{split}
\left[ \intO{ B(\vr) } \right]_{t = 0}^{t = \tau} 
= &- \int_{0}^{\tau} \intO{ \Div (\vr \vu) B'(\vr) } \dt - \ep \int_{0}^{\tau} \intO{ |\Grad \vr|^2  B''(\vr) } \dt  
\\ & + \int_{0}^{\tau} \int_{\partial \Omega} B'(\vr) (\vr - \vr_B) [\vu_B \cdot \vc{n}]^- \ {\rm d} S_x \ \dt
\end{split}
\end{equation}
holds for any $0 \leq \tau \leq T$. As consequence we obtain strict positivity of $\vr$ on condition that $\vr_B$, $\vr_0$ enjoy the same property,
\begin{equation}\label{min}
{\rm ess} \inf_{t,x} \vr(t,x) \geq  
\min \left\{ \min_\Omega \vr_0 ; \min_{\Gamma_{\rm in}} \vr_B \right\} 
C(\| \vu_B \|_{L^\infty((0,T) \times \Omega)} ) \exp \left( -T \| \Div \vu \|_{L^\infty((0,\tau) \times \Omega)} \right).
\end{equation}
Finally, since $\vr_B$, $\vu_B$ are smooth and time independent vector fields, we have the estimates
\begin{equation} \label{E4a}
\| \partial_t \vr \|_{L^2((0,T) \times \Omega)} + 
\ep {\rm ess} \ \sup_{\tau \in (0,T)} \left\| \Grad \vr (\tau, \cdot) \right\|^2_{L^2(\Omega)} \leq c,
\end{equation} 
with the constant depending only on the data, for details see \cite{AbFeNo}.

\subsection{Regularization of the pressure}
The next level of approximation consists in replacing the pressure in the momentum equation \eqref{mom-e} by an isoentropic pressure 
\begin{equation}\label{pi-m}
\pi_m(\vr)=\vr^{\am}. \end{equation}
where  $\am$ is  a sequence of nonnegative real numbers such that $\am\to + \infty$.

\subsection{Galerkin approximations of the momentum balance}
Finally, we look for approximate velocity field in the form 
\[
\vu = \vv + \vu_B, \ \vv \in C([0,T]; X_n). 
\]
such that the approximated momentum balance in the weak form is given by
\begin{equation} \label{E5}
\begin{split}
&\left[ \intO{ \vr \vu \cdot \bfphi } \right]_{t=0}^{t = \tau} \\&= 
\int_0^\tau \intO{ \Big[ \vr \vu \cdot \partial_t \bfphi + \vr \vu \otimes \vu : \Grad \bfphi 
+ \pi_m(\vr) \Div \bfphi - \partial F_\delta (\Ds \vu)  : \Grad \bfphi \Big] }
\\&- \ep \int_0^\tau \intO{ \Grad \vr \cdot \Grad \vu \cdot \bfphi } \dt \mbox{ for any } \bfphi \in C^1([0,T]; X_n)
\end{split}
\end{equation}
and with the initial condition
\begin{equation} \label{E6}
\vr \vu(0, \cdot) = \vr_0 \vu_0, \ \vu_0 = \vv_0 + \vu_B, \ \vv_0 \in X_n.
\end{equation}

For fixed parameter $n, m\in \mathbb{N}$, $\delta > 0$, $\ep > 0$, the  approximating system is given by  the parabolic problem \eqref{E1}--\eqref{E3}, the pressure \eqref{pi-m} and  the Galerkin approximation \eqref{E5}, \eqref{E6}. 
The existence of the approximate solutions  can be proved through standard arguments. 
For $\vu = \vu_B + \vv$, $\vv \in C([0,T]; X_n)$, we identify the unique solution $\vr = 
\vr [\vu]$ of \eqref{E1}--\eqref{E3} and plug it as $\vr$ in \eqref{E5}. The unique solution $\vu = \vu[\vr]$ of \eqref{E5} defines a mapping 
\[
\mathcal{T}: \vv \in C([0,T]; X_n) \mapsto 
\mathcal{T}[\vc{v}] = (\vu[\vr] - \vu_B)  
\in C([0,T]; X_n).
\]
The first level approximate solutions $\vr = \vr_{\delta, \ep, m, n}$, $\vu = \vu_{\delta, \ep,m,  n}$ are obtained via a fixed point 
argument to the mapping $\mathcal{T}$. 

Moreover, the approximate energy inequality 
\begin{equation} \label{E7}
\begin{split}
&\left[ \intO{\left[ \frac{1}{2} \vr |\vu - \vu_B|^2 +\Pi_m(\vr) \right] } \right]_{t = 0}^{ t = \tau} + 
\int_0^\tau \intO{ \partial F_\delta( \Ds \vu) : \Ds \vu } \dt \\ 
&+\int_0^\tau \int_{\Gamma_{\rm out}}  \Pi_m(\vr)  \vu_B \cdot \vc{n} \ \D S_x \dt\\
&- \int_0^\tau \int_{\Gamma_{\rm in}} \left[ \Pi_m(\vr_B) - \Pi_m'(\vr) (\vr_B - \vr) - \Pi_m(\vr) \right] \vu_B \cdot \vc{n} 
\ {\rm d}S_x \ \dt + \ep \int_0^\tau \intO{ \Pi_m''(\vr) |\Grad \vr|^2 } \dt\\
\leq
&- 
\int_0^\tau \intO{ \left[ \vr \vu \otimes \vu + \pi_m(\vr)  \mathbb{I} \right]  :  \Grad \vu_B } \dt + \int_0^\tau \intO{ {\vr} \vu  \cdot \vu_B \cdot \Grad \vu_B  } 
\dt\\ &+ \int_0^\tau \intO{ \partial F_\delta( \Ds \vu) : \Grad \vu_B } \dt - \int_0^\tau \int_{\Gamma_{\rm in}}  \Pi_m(\vr_B) \vu_B \cdot \vc{n} \ \D S_x \dt  
\end{split}
\end{equation}
holds for any $0 \leq \tau \leq T$ and for any $n, m\in \mathbb{N}$, $\delta > 0$, $\ep > 0$ and where 
$$\Pi_m(\vr) = \vr \int_0^\vr \frac{\pi_m(z)}{z^2}{\rm d}z= \frac{1}{\am -1}\,\vr^\am.$$

\subsection{Approximated initial and boundary data}
By using for simplicity only the subscript $n$ to denote the sequences of approximating initial and boundary data we assume that
\begin{equation}
\vr_n|_{t=0} = \vr_{n_0}, \quad \vr_n \vu_n|_{t=0} =\vm_{n_0}, 
\label{i1}
\end{equation}
where
\begin{eqnarray}
& & 0 \le \vr_{n_0} \quad \mbox{a.e}, \quad  \vr_{n_{0}} \in L^1(\Omega) \cap L^{\alpha_m}(\Omega),  
\quad \int (\vr_{n_0})^{\alpha_m} dx \le c \alpha_m \,\, \mbox{for some}\,\, c, \label{i5}\\
& & \quad \quad \quad  \vm_{n_0} \in L^{\frac{2 \alpha_{m}}{\alpha_m + 1}} (\Omega), \quad \vu_{n_0} = \frac{\vm_{n_0}}{\vr_{n_0}}\,\, \mbox{on} \,\, \{\vr_{n_0} >0\}, \quad \vu_{n_0} = 0\,\, \mbox{on} \,\, \{\vr_{n_0} = 0\} \nonumber \\
&& \vr_{n_0} |\vu_{n_0}|^2  \,\, \mbox{is bounded in}\,\,  L^1(\Omega).
\end{eqnarray}
Furthermore we assume that
\begin{equation}
M_{m}=\int_{\Omega}\!\!\!\!\!\!\!-\vr_{n_0},\quad  0<M_{m}<M<1,\quad M_{m}\to M.
\label{i4}
\end{equation}
\begin{equation}\nonumber
\vr_{n_0}\vu_{n}\rightharpoonup \vm_{0}\quad \text{weakly in $L^{2}(\Omega)$,}
\end{equation}
\begin{equation}
\vr_{n_0}\rightharpoonup \vr_{0}\quad \text{weakly in $L^{1}(\Omega)$.}
\label{i3}
\end{equation}
For the boundary data we assume the same conditions as in  Section \ref{boundary data},
\begin{align} \label{P5bis}
\vu_{n}|_{\partial \Omega} &= \vu_B,\\
\vu_B&\in C^2(\partial\Omega;\mathbb{R}^d), \quad  \int_{\partial\Omega} \vu_B\cdot\vn \geq 0, \label{bdu-bis}\\
 \vr_{n}|_{\Gamma_{in}} &= \vr_B, \quad  \vr_B \in C^1(\partial \Omega), \quad  0<\vr_B\leq 1  \mbox{ on } \Gamma_{\rm in}, \quad \inf_{\Gamma_{\rm in}} \vr_B>0. \label{hyp-rb-bis}
\end{align}

\section{Existence of dissipative solutions}
\label{existence}
This section is devoted to the proof of the Theorem \ref{T1}. Therefore we have to perform different limits in the approximating system introduced in the previous Section \ref{E}.  We perform first the limit as $\delta\to 0$, then $\ep\to 0$, $m\to \infty$ and at the end the Galerkin limit. This choice is due to the fact that we need to take advantage on the additional estimates that comes from the finite dimensional space. At this level we assume that the initial and boundary data have enough regularity.

\subsection{Limit as $\delta \to 0$}
We take the limit as $\delta \to 0$ in the regularization of the potential $F_{\delta}$ keeping the remaining parameters fixed. 
From  the energy inequality \eqref{E7} and the property \eqref{F1a}
we deduce a uniform bound for $\delta \to 0$ on the traceless part of the symmetric velocity gradient, 
\[
\left\| \Ds \vu - \frac{1}{d} \Div \vu \mathbb{I} \right\|_{L^q((0,T) \times \Omega; \mathbb{R}^{d \times d})} \leq c,\ q > 1. 
\]
 Employing a generalized Korn inequality we obtain a bound uniform in $\delta$ for
$\left\|  \Grad \vu  \right\|_{L^q((0,T) \times \Omega; \mathbb{R}^{d \times d})},$
then using also the standard Poincar\' e inequality (recalling $\vu=\vu_B$ on $\partial\Omega$), yields 
\begin{equation} \label{E10}
\left\| \vu \right\|_{L^q(0,T; W^{1,q}(\Omega; \mathbb{R}^d))} \leq C \quad  \mbox{for some}\ q > 1.
\end{equation}
Since $n$ is fixed all norms are equivalent on the finite--dimensional space $X_n$. In particular, 
\[
\left\| \Grad \vu \right\|_{L^\infty((0,T) \times \Omega; \mathbb{R}^{d \times d})}\leq C \ \mbox{ uniformly for } \delta \searrow 0.
\]
Hence it is standard to perform the limit $\delta \to 0$ and deduce from  \eqref{E5} that it holds
\begin{equation} \label{E11}
\begin{split}
&\left[ \intO{ \vr \vu \cdot \bfphi } \right]_{t=0}^{t = \tau} = 
\int_0^\tau \intO{ \Big[ \vr \vu \cdot \partial_t \bfphi + \vr \vu \otimes \vu : \Grad \bfphi 
+ \pi_m(\vr) \Div \bfphi - \mathbb{S}  : \Grad \bfphi \Big] }\\
&- \ep \int_0^\tau \intO{ \Grad \vr \cdot \Grad \vu \cdot \bfphi } \dt,\ 
\mathbb{S} \in L^\infty((0,T) \times \Omega; R^{d \times d}_{\rm sym}), 
\end{split}
\end{equation}
for any  $\bfphi \in C^1([0,T]; X_n)$ and the energy inequality in the form 
\begin{equation} \label{E12}
\begin{split}
&\left[ \intO{\left[ \frac{1}{2} \vr |\vu - \vu_B|^2 +  \Pi_m(\vr)\right] } \right]_{t = 0}^{ t = \tau} + 
\int_0^\tau \intO{\left[ F(\Ds \vu) + F^*(\mathbb{S}) \right] } \dt \\ 
&+\int_0^\tau \int_{\Gamma_{\rm out}} \Pi_m(\vr) \vu_B \cdot \vc{n} \ \D S_x \dt\\
&- \int_0^\tau \int_{\Gamma_{\rm in}} \left[ \Pi_m(\vr_B) - \Pi_m'(\vr) (\vr_B - \vr) - \Pi_m(\vr) \right] \vu_B \cdot \vc{n} 
\ {\rm d}S_x \ \dt + \ep \int_0^\tau \intO{ \Pi_m''(\vr) |\Grad \vr|^2 } \dt
\\
\leq
&- 
\int_0^\tau \intO{ \left[ \vr \vu \otimes \vu + \pi_m(\vr)  \mathbb{I} \right]  :  \Grad \vu_B } \dt + \int_0^\tau \intO{ {\vr} \vu  \cdot \vu_B \cdot \Grad \vu_B  } 
\dt\\ &+ \int_0^\tau \intO{ \mathbb{S} : \Grad \vu_B } \dt - \int_0^\tau \int_{\Gamma_{\rm in}} \Pi_m(\vr_B)  \vu_B \cdot \vc{n} \ \D S_x \dt. 
\end{split}
\end{equation}

\subsection{Vanishing viscosity limit: $\ep \to 0$}
The purpose of this section is to prove the limit as $\ep \to 0$. One of the most delicate point is to recover the a.e. convergence of the density in order to perform the limit of the pressure $\pi_{m}$.

Since $n$ is still fixed, we have 
\begin{equation} \label{E13}
\| \vu_\ep \|_{L^\infty(0,T; W^{1,\infty}(\Omega))} \leq C, 
\end{equation}
yielding, as in \eqref{max} and \eqref{min}, the uniform bounds on the density 
\begin{equation} \label{E14}
0 < \underline{\vr} \leq \vr_\ep(t,x) \leq \Ov{\vr} \ \mbox{for all}\ (t,x) \in [0,T] \times \Ov{\Omega}.
\end{equation} 
Since at this stage $\vr_{\ep}$ possess a well defined trace on $\partial \Omega$, for fixed $n, m > 0$, by virtue of the energy inequality \eqref{E12}, we have
\begin{equation} \label{E15a}
\sqrt{\ep} \|\Grad \vr_\ep \|_{L^2((0,T) \times \Omega)}\leq C, 
\end{equation}
and 
\begin{equation} \label{E15}
\sup_{\tau \in [0,T]} \left\| \vu_\ep(\tau, \cdot) \right\|_{W^{1, \infty}(\Omega; \mathbb{R}^d)} \leq C.
\end{equation}
As consequence 
\begin{equation} \label{dens}
\begin{split}
\vre \to \vr \ &\mbox{weakly-(*) in}\ L^\infty((0,T) \times \Omega)\\ 
\ &\mbox{weakly in}\ C_{\rm weak}([0,T]; L^s(\Omega)) \ \mbox{for any}\ 1 < s < \infty, \\
\ &\mbox{strongly in}\  L^{2}((0,T); H^{-1}(\Omega)), 
\end{split}
\end{equation}
passing to a suitable subsequence. The second convergence follows form the weak bound on the 
time derivative $\partial_t \vre$ obtained from equation \eqref{E4}. Moreover have 
\[
\vre \to \vr \ \mbox{weakly-(*) in}\ L^\infty((0,T) \times \partial \Omega; \D S_x)
\]
and the limit density admits the same upper and lower bounds as in \eqref{E14}. 

Similarly, 
\begin{equation} \label{E15b}
\vue \to \vu \ \mbox{weakly-(*) in}\ L^\infty(0,T; W^{1,\infty}(\Omega; \mathbb{R}^d)), 
\end{equation}
and 
\begin{equation} \label{rho-u}
\vre \vue \to \vc{m} \ \mbox{weakly-(*) in} \ L^\infty((0,T) \times \Omega; \mathbb{R}^d)).
\end{equation}
An abstract version of Arzela--Ascoli theorem yields 
\begin{equation} \label{E15c}
\vm = \vr \vu \ \mbox{a.a. in}\ (0,T) \times \Omega.
\end{equation}

\subsubsection{The limit in the approximate equation of continuity}

By virtue of the established convergences we can take the limit in the regularized equation of continuity \eqref{E4}:
\begin{equation} \label{E16}
\begin{split}
\left[ \intO{ \vr \varphi } \right]_{t = 0}^{t = \tau}&= 
\int_0^\tau \intO{ \Big[ \vr \partial_t \varphi + \vr \vu \cdot \Grad \varphi  \Big] }
\dt \\ 
&- \int_0^\tau \int_{\Gamma_{\rm out}} \varphi \vr \vu_B \cdot \vc{n} \ {\rm d} S_x \dt - 
\int_0^\tau \int_{\Gamma_{\rm in}} \varphi \vr_B \vu_B \cdot \vc{n}  \ {\rm d}S_x \dt,\ 
\vr(0, \cdot) = \vr_0 
\end{split}
\end{equation}
for any $\varphi \in C^1([0,T] \times \Ov{\Omega})$. We point out that the quantity 
\[
\vm = \left\{ \begin{array}{l} \vr \vu_B \cdot \vc{n} \ \mbox{on}\ \Gamma_{\rm out},\\ \\
\vr_B \vu_B \cdot \vc{n} \ \mbox{on}\ \Gamma_{\rm in} \end{array} \right.  
\]
is the normal trace of the divergenceless vector field $[\vr, \vr \vu]$ on the lateral boundary 
$(0,T) \times \partial \Omega$ in the sense of Chen-Torres-Ziemer \cite{ChToZi}.

\subsubsection{The limit in the approximate momentum equation}

By the assumption \eqref{S3},  $F$ is proper convex with ${\rm Dom}[F] = R^{d \times d}_{\rm sym}$ hence it  holds that $F^*$ is a superlinear function. From the energy inequality \eqref{E12} we deduce that
 $$\intTO{F^*(\mathbb{S_\ep})}$$  is uniformly bounded as $\ep\to 0$, thus
\[
\mathbb{S}_\ep \to \mathbb{S} \ \mbox{weakly in}\ L^1((0,T) \times \Omega; R^{d \times d}). 
\]
Next, from \eqref{E11} it follows that 
\[
\partial_t 
P_n[ \vre \vue ] \ \mbox{bounded in}\ L^2(0,T; X_n), 
\]
where $P_n: L^2 \to X_n$ is the associated orthogonal projection. By using \eqref{E15b}, \eqref{E15c}, 
we may infer that 
\[
\vre \vue \otimes \vue \to \vr \vu \otimes \vu \ \mbox{weakly-(*) in}\ 
L^\infty((0,T) \times \Omega). 
\]
Now  we may let $\ep \to 0$ in \eqref{E11} obtaining 
\begin{equation} \label{E20}
\begin{split}
\left[ \intO{ \vr \vu \cdot \bfphi } \right]_{t=0}^{t = \tau} &= 
\int_0^\tau \intO{ \Big[ \vr \vu \cdot \partial_t \bfphi + \vr \vu \otimes \vu : \Grad \bfphi 
+ {\overline{\pi_m(\vr)}} \Div \bfphi - \mathbb{S}  : \Grad \bfphi \Big] }
\end{split}
\end{equation}
for any $\bfphi \in C^1([0,T]; X_n)$, where  ${\overline{\pi_m(\vr)}}\in L^\infty((0,T) \times \Omega)$ denotes the weak limit as $\ep \to 0$ of the sequence $\pi_{m}(\vr_{\ep})=\vr_{\ep}^{\alpha_{m}}$. The next step is to  show that 
$${\overline{\pi_m(\vr)}}= \pi_m(\vr) \in L^\infty((0,T) \times \Omega).$$
In order to motivate it, we suitably modify the standard arguments by Feireisl and the new setting of in/out-flow boundary condition by Chang-Jin-Novotn\'{y} \cite{ChJiNo}. (Note that in \cite{ChJiNo} the limit in the viscous approximation of the continuity equation is performed after the limit in the Galerkin approximation, while the situation here is different and we still have the advanteges of keeping the finite dimension $n$ fixed.)

Since $n$ is fixed it holds that
$$ \Div \vu_\ep \to \Div \vu  \ \mbox{weakly in}\ 
L^\infty((0,T), L^2(\Omega)),
$$
then using \eqref{dens} we can infer that
\begin{equation}\label{rho-div-e}
\vr_\ep\Div \vu_\ep \to \vr\Div \vu  \ \mbox{weakly-(*) in}\ 
L^\infty((0,T) \times \Omega).
\end{equation}
Since $\vr\Div \vu$ is integrable we get from the general theory of transport equations developed by DiPerna-Lions that the limit $[\vr, \vu]$ is a renormalized solution in the sense specified in \eqref{renorm}. In particular the choice $B: z\mapsto z\log z$ leads to 
\begin{equation}\label{renorm3}
\begin{split}
\left[ \intO{{\vr\log\vr} } \right]_{t = 0}^{t = \tau} 
= &- \int_{0}^{\tau} \intO{ \Div (\vr \vu) B'(\vr) } \dt  
\\ & + \int_{0}^{\tau} \int_{\partial \Omega} B'(\vr) (\vr - \vr_B) [\vu_B \cdot \vc{n}]^- \ {\rm d} S_x \ \dt,
\end{split}
\end{equation}

On the other hand, let $B$ be a smooth convex function and let us consider \eqref{renorm} for any $\ep>0$ fixed, using that $ \ep \int_{0}^{\tau} \intO{ |\Grad \vr_\ep|^2  B''(\vr_\ep) } \dt  $ is non-negative it yields that
\begin{equation}\label{renorm2}
\begin{split}
\left[ \intO{ B(\vr_\ep) } \right]_{t = 0}^{t = \tau} 
\leq  &- \int_{0}^{\tau} \intO{ \Div (\vr_\ep \vu_\ep) B'(\vr_\ep) } \dt
\\ & + \int_{0}^{\tau} \int_{\partial \Omega} B'(\vr_\ep) (\vr_\ep - \vr_B) [\vu_B \cdot \vc{n}]^- \ {\rm d} S_x \ \dt,
\end{split}
\end{equation}
then choosing $B$ as the function $z\mapsto z\log z$ and taking the limit as $\ep \to 0$ we deduce 
\begin{equation}\label{renorm4}
\begin{split}
\left[ \intO{ \overline{\vr\log\vr} } \right]_{t = 0}^{t = \tau} 
\leq &- \int_{0}^{\tau} \intO{ \Div (\vr \vu) B'(\vr) } \dt 
\\ & + \int_{0}^{\tau} \int_{\partial \Omega} B'(\vr) (\vr - \vr_B) [\vu_B \cdot \vc{n}]^- \ {\rm d} S_x \ \dt.
\end{split}
\end{equation}
Finally subtracting \eqref{renorm4} and \eqref{renorm3} we have
\begin{equation}\label{renorm5}
\begin{split}
\left[ \intO{ \overline{\vr\log\vr} - \vr\log\vr } \right]_{t = 0}^{t = \tau} 
\leq0,
\end{split}
\end{equation}
but, since the function  $B: z\mapsto z\log z$ is convex it holds $\overline{\vr\log\vr} - \vr\log\vr\geq 0 $ a.e. in $(0,T)\times\Omega$.
Therefore
$$ \overline{\vr\log\vr} = \vr\log\vr \mbox{ a.e. in } (0,T)\times\Omega,$$
and 
\begin{equation}
\vr_\ep\to \vr \ \mbox{ a.e. in } (0,T)\times\Omega.
\end{equation}
Consequently  for any $m, n\in \mathbb{N}$ fixed,  \eqref{E20} reads as 
\begin{equation} \label{MM}
\begin{split}
\left[ \intO{ \vr \vu \cdot \bfphi } \right]_{t=0}^{t = \tau} &= 
\int_0^\tau \intO{ \Big[ \vr \vu \cdot \partial_t \bfphi + \vr \vu \otimes \vu : \Grad \bfphi 
+ {\pi_m(\vr)} \Div \bfphi - \mathbb{S}  : \Grad \bfphi \Big] }
\end{split}
\end{equation}
for any $\bfphi \in C^1([0,T]; X_n)$.

Finally, by using the weak lower semicontinuity of the potentials $F$ and $F^{\ast}$ we perform the limit in the energy \eqref{E12} and we obtain  the energy inequality holds in the form
\begin{equation} \label{En-ep}
\begin{split}
&\left[ \intO{\left[ \frac{1}{2} \vr |\vu - \vu_B|^2 + \Pi_m(\vr)\right] } \right]_{t = 0}^{ t = \tau} + 
\int_0^\tau \intO{\left[ F(\Ds \vu) + F^*(\mathbb{S}) \right] } \dt \\ 
&+\int_0^\tau \int_{\Gamma_{\rm out}} \Pi_m(\vr) \vu_B \cdot \vc{n} \ \D S_x \dt\\
&- \int_0^\tau \int_{\Gamma_{\rm in}} \left[ \Pi_m(\vr_B) - \Pi_m'(\vr) (\vr_B - \vr) - \Pi_m(\vr) \right] \vu_B \cdot \vc{n} 
\ {\rm d}S_x \ \dt 
\\
\leq
&- 
\int_0^\tau \intO{ \left[ \vr \vu \otimes \vu + \pi_m(\vr)  \mathbb{I} \right]  :  \Grad \vu_B } \dt + \int_0^\tau \intO{ {\vr} \vu  \cdot \vu_B \cdot \Grad \vu_B  } 
\dt\\ &+ \int_0^\tau \intO{ \mathbb{S} : \Grad \vu_B } \dt - \int_0^\tau \int_{\Gamma_{\rm in}} \Pi_m(\vr_B)  \vu_B \cdot \vc{n} \ \D S_x \dt. 
\end{split}
\end{equation}

\subsection{Blow up of the adiabatic exponent: the limit as $m\to + \infty$}
In this section we perform the limit in the adiabatic coefficient $\alpha_{m}$. At this level, since $m\to + \infty$ we lose the uniform bounds on the pressure term $\pi_{m}(\vr_{\ep})=\vr_{\ep}^{\alpha_{m}}$. Hence one of the main issue of this section is to recover the uniform integrability bounds of the density. We start by collecting the uniform bounds that follow directly from the energy inequality \eqref{En-ep}. Employing the growth condition on the potential $F$ and recalling the de la Vall\'{e}-Poussin criterion as $F^*$ is a superlinear function, we get that
\begin{align}
&\vu_m\to \vu  \mbox{ weakly in }\ L^q(0,T; W^{1,q}(\Omega; \mathbb{R}^d)), \label{u-m-infty}\\
&\mathbb{S}_m\to \mathbb{S}  \mbox{ weakly in }\ L^1((0,T)\times \Omega; \mathbb{R}^{d\times d}).
\end{align}
Since the finite dimension $n$ of  the Galerkin level approximation is still fixed, $\vu_m$ enjoys that
\begin{equation}\label{u-m}
\vu_m \to \vu \ \mbox{weakly-(*) in}\ L^\infty(0,T; W^{1,\infty}(\Omega; \mathbb{R}^d)).
\end{equation}
Consequently, after taking the limit as $\ep \to 0$ in \eqref{max}, we deduce that for $\vr_{m}$ it is still valid the following bound
\begin{equation}\begin{split}
\|\vr_m\|_{L^\infty((0,\tau) \times \Omega)} \leq &
\max \left\{ \| \vr_0 \|_{L^\infty(\Omega)}; \| \vr_B \|_{L^\infty((0,T) \times \Gamma_{in})};
\| \vu_B \|_{L^\infty((0,T) \times \Omega)} \right\}\\
&\cdot \exp \left( \tau \| \Div \vu \|_{L^\infty((0,\tau) \times \Omega)} \right),
\end{split}
\label{dens-m1}
\end{equation}
and it holds (up to a subsequence) that
\begin{equation} \label{dens-m}
\vr_m\to \vr \ \mbox{weakly-(*) in}\ L^\infty((0,T) \times \Omega) 
\ \mbox{and weakly in}\ C_{\rm weak}([0,T]; L^s(\Omega)) \ \mbox{for any}\ 1 < s < \infty.
\end{equation}
Moreover,
 \begin{equation} \label{rho-u-m}
\vr_m \vu_m \to \vr\vu \ \mbox{weakly-(*) in} \ L^\infty((0,T) \times \Omega; \mathbb{R}^d)).
\end{equation}
The limit $[\vr, \vu]$ enjoys the weak formulation of the continuity equation in the form \eqref{E16} and also it holds the total mass inequality
\begin{equation}\label{total-mass}
 \intO{\vr(t)} + \int_{\Gamma_{\rm out}} \vr  \vu_B\cdot \vn \,\D S_x \dt  \leq \intO{\vr_0}- \int_0^t \int_{\Gamma_{\rm in}} \vr_B \vu_B\cdot \vn \,\D S_x \dt \mbox{ for any } t\in [0, T].
\end{equation}

\subsubsection{Boundedness of the density} 
Recalling the definition of the pressure $\pi_m$ in \eqref{pi-m} and computing the pressure potential $\Pi_m$,  we deduce from the energy inequality \eqref{En-ep} that
\begin{equation} \label{Press-pot}
\intO{\vr_m(t)^{\am}}\leq C \am. \end{equation}
Unfortunately this estimate is not uniform in $m$. Our aim, now, is to prove that $\vr_m^{\am}$ is uniformly bounded and that $(\vr_m -1)_+\to 0 $ in $L^\infty(0, T; L^p(\Omega))$ for any $1\leq p<+\infty$. 

For any $1<p<\infty$ and for any $m$ such that $\am>p$ using \eqref{Press-pot} we have 
\begin{equation}\label{interpolation}
\|\vr_m\|_{L^\infty(0,T; L^p(\Omega))}\leq \|\vr_m\|_{L^\infty(0,T; L^1(\Omega))}^{\vartheta_m} \|\vr_m\|_{L^\infty(0,T; L^{\am}(\Omega))}^{1-\vartheta_m}\leq M_m^{\vartheta_m} (C\am)^{\frac{1-\vartheta_m}{\am}}
\end{equation}
where $$\frac{1}{p}=\vartheta_m + \frac{1-\vartheta_m}{\am}$$
and 
$ M_m $ is defined in \eqref{i4}.
Since $\alpha_m\to +\infty$ as $m\to +\infty$ we deduce that $\vartheta_m\to \frac{1}{p}$ and 
$$\|\vr\|_{L^\infty(0,T; L^p(\Omega))}\leq \liminf_{m\to +\infty}\|\vr_m\|_{L^\infty(0,T; L^p(\Omega))} \leq M^{\frac{1}{p}}.$$
Next, letting $p\to \infty$ we obtain
\begin{equation}\label{limit-p}
\|\vr\|_{L^\infty((0,T) \times\Omega)}\leq \liminf_{p\to +\infty}\|\vr\|_{L^\infty(0,T; L^p(\Omega))}\leq 1.
\end{equation}
Now we introduce $\phi_m =(\vr_m -1)_+$. A direct application of \eqref{Press-pot} implies that
$$ \intO{(1+\phi_m)^{\am} 1_{\phi_m >0}} \leq \intO{\vr_m^{\am}}\leq C \am.$$
By virtue of the algebraic inequality 
$$ (1+ x)^k\geq 1 + C k^p x^p \mbox{ for any } x\geq 0$$
valid for any $p>1$ and $k=k(p)$ large enough, with $C=C(p)$ positive constant, it yields that
$$ \intO{\phi_m^p}\leq \frac{C}{\am^{p-1}},$$
which ensures that
\begin{equation}\label{positivepart}
(\vr_m -1)_+ \to 0 \mbox{ as } m\to +\infty \quad \text{in $L^{\infty}((0,T), L^{p}(\Omega))$, $1\leq p<+\infty$.}
\end{equation}
Therefore the limit $\vr$ is bounded as follows
\begin{equation}\label{rho-limitata}
0\leq \vr\leq 1 \mbox{ a.e. in } \Omega.
\end{equation}
Let us observe that from \eqref{En-ep} we also deduce that 
\begin{equation} \label{Press-pot-boundary}
\int_0^T \int_{\Gamma_{\rm out}} \vr_m(t)^{\am}\vu_B \cdot \vc{n} \ \D S_x \dt \leq C \am
 \end{equation}
with $C$ a positive constant depending on the data. Also from \eqref{total-mass} $\vr$ is uniformly bounded in $L^\infty(0,T; L^1(\Gamma_{\rm out}, |\vu_B\cdot\vn|  \D S_x ))$. Consequently, we may repeat the interpolation and the limit procedure performed above for $\vr_m$ on $\Omega$ (see \eqref{interpolation}--\eqref{limit-p}) and deduce that
\begin{equation}
\|\vr\|_{L^\infty((0,T)\times \Gamma_{\rm out}, |\vu_B\cdot\vn|  \D S_x)}\leq 1,
\end{equation}
and
\begin{equation}\label{rho-bordo-limitata}
0\leq \vr\leq 1 \mbox{ a.e. in } (0,T)\times\Gamma_{\rm out}
\end{equation}
with respect to the measure on the boundary $|\vu_B\cdot\vn|  \D S_x$.

\subsubsection{Uniform integrability of the pressure}

Let us consider the quantity 
$$\bfphi(t,x)= \eta(t)\, \phi(x) \Grad \Del^{-1}(\phi(x) \vr_m(t,x)), \ \eta\in C^\infty_{\rm c}(0,T), \  \phi\in C^\infty_{\rm c}(\Omega)$$
as test function in the approximate momentum equation \eqref{MM} and where $\Del^{-1}$ is the inverse of the Laplace operator on the space $\mathbb{R}^d$ with $d=2,3$ defined via the convolution with the Poisson kernel. 
We obtain
\begin{equation} \label{MM2}
\int_0^\tau\eta \intO{ \phi^2\,\pi_m(\vr_m)\vr_m}=\sum_{j=1}^6 I_j,
\end{equation}
where
\begin{align*}
&I_1=-\int_0^\tau\eta  \intO{\pi_m(\vr_m)\Grad\phi\cdot  \Grad \Del^{-1}(\phi \vr_m)},\\
&I_2= -\int_0^\tau\eta \intO{\phi\,\vr_m \vu_m \cdot \partial_t(\Grad \Del^{-1}(\phi \vr_m)) }\dt,\\
&I_3= - \int_0^\tau\eta \intO{\vr_m \vu_m \otimes \vu_m \cdot \Grad \phi \cdot \Grad \Del^{-1}(\phi \vr_m) }\dt,\\
&I_4= - \int_0^\tau\eta \intO{\phi\, \vr_m \vu_m \otimes \vu_m : \Grad^2\Del^{-1}(\phi\vr_m) }\dt,\\
&I_5= \int_0^\tau\eta \intO{ \mathbb{S}_m  \cdot \Grad \phi \cdot \Grad \Del^{-1}(\phi \vr_m)}\dt,\\
&I_6= \int_0^\tau\eta \intO{ \phi\,\mathbb{S}_m  : \Grad^2 \Del^{-1}(\phi \vr_m)}\dt.
\end{align*}
Our aim is to show that all the integrals on the right hand side of \eqref{MM2} are bounded in terms of the uniform bounds established in \eqref{En-ep} -- \eqref{rho-u-m} and \eqref{positivepart}. We need to consider the time derivative
$$ \partial_t(\Grad \Del^{-1}(\phi(x) \vr_m(t,x)))= \Grad \Del^{-1}(\phi(x) \partial_t\vr_m(t,x)).$$
Since $\vr_m$ satisfies the approximate continuity equation \eqref{E16} we have 
$$\phi\partial_t\vr_m = -\phi \Div(\vr_m\vu_m) \mbox{ in } \mathcal{D}'((0,T)\times\Omega).$$
Being $\phi$ compactly supported the boundary conditions become irrelevant.
We deduce that
\begin{equation}
\Grad \Del^{-1}(\phi \partial_t\vr_m)=-\Grad\Del^{-1}\Div(\phi\vr_m\vu_m) + \Grad\Del^{-1}[\vr_m\vu_m\cdot\Grad\phi].
\end{equation}
In view of  \eqref{rho-u-m} and of the properties of the operator $\Del^{-1}$  we get
\begin{equation}
{\rm ess}\sup_{(0,T)}\|-\Grad\Del^{-1}\Div(\phi\vr_m\vu_m) + \Grad\Del^{-1}[\vr_m\vu_m\cdot\Grad\phi]\|_{L^\infty(\Omega; \mathbb{R}^d)} \leq C(\phi){\rm ess}\sup_{(0,T)}\|\vr_m\vu_m\|_{L^\infty(\Omega; \mathbb{R}^d)}.
\end{equation}
Similarly, employing the uniform bounds \eqref{dens-m1}, \eqref{dens-m} we observe that 
\begin{equation}
\|\Grad \Del^{-1}(\phi \vr_m)\|_{L^\infty((0,T)\times\Omega; \mathbb{R}^d)}\leq C(\phi),
\end{equation}
and 
\begin{equation}
{\rm ess}\sup_{(0,T)}\|\Grad^2\Del^{-1}(\phi\vr_m)\|_{L^s(\Omega)}\leq C(\phi){\rm ess}\sup_{(0,T)}\|\vr_m\|_{L^s(\Omega)}
\end{equation}
Going back in \eqref{MM2} we obtain that for any compact set $K\subset\Omega$ it holds
\begin{equation}\label{pi-m-equi-K}
\|\eta\pi_m(\vr_m)\vr_m\|_{L^1((0,T)\times K)}\leq C({\rm data}, n, K). 
\end{equation}

Hence, for fixed $n\in\mathbb{N}$ the sequence $\pi_m(\vr_m)$ is equintegrable in $L^1(Q)$ for any $Q\subset\subset (0,T)\times\Omega$. 
 
 Now, in order to extend \eqref{pi-m-equi-K} up to the boundary we use a similar strategy with a test function
 $$\bfphi= \mathcal{B}[{\Phi}], \ \Phi\in L^q(\Omega), \ \intO{\Phi}=0,$$
 where $\mathcal{B}$ is  the Bogovskii operator that acts as the inverse of the divergence operator, for its main properties see Section 0.8 in \cite{FNopen}.
 
  \begin{equation}\label{Bog-boundary}
\begin{split}
\intTO{&\pi_m(\vr_m) \Phi}= \left[\intO{\vr_m\vu_m\cdot \mathcal{B}[{\Phi}]}\right]_{t=0}^{t=T}
\\&-\intTO{ \vr_m \vu_m \otimes \vu_m : \Grad  \mathcal{B}[{\Phi}]}
+ \intTO{\mathbb{S}_m  : \Grad  \mathcal{B}[{\Phi}]}.
\end{split}
\end{equation}
 By virtue of the properties of the Bogovskii operator and of the uniform boundedness of the sequences $\vr_m$ and $\vr_m\vu_m$ in the spaces $L^p((0,T)\times\Omega)$ for any $p$ we are in position to claim that the integral on the left-hand side of \eqref{Bog-boundary} is
bounded. Choosing a suitable $\Phi$ singular near the boundary, we  combine \eqref{Bog-boundary} with \eqref{pi-m-equi-K} to conclude 
  \begin{equation}\label{Bog-boundary2}\intTO{\pi_m(\vr_m) {\rm dist}^{-\omega}[x, \partial\Omega]}\leq C \mbox{ for some } \omega>0, \mbox{ for any } m\in\mathbb{N}
  \end{equation}
 then \eqref{pi-m-equi-K}  and \eqref{Bog-boundary2} ensure
 the equintegrability of the pressure sequence $\pi_m(\vr_m)$. 
We established that
\begin{equation}\label{vr-m+1}
\int_{Q}{\vr_m^{\am+1}}\leq C({\rm data}, n)|Q|, \mbox{ for any }Q\subset\subset (0,T)\times\Omega.
\end{equation}
 Let us prove that $\vr_m^\am$ is equintegrable in $L^1$.  We have  that
\begin{equation}
\int_{J\times K} \vr_m^\am \dx\dt \leq \left(\int_{J\times K}\vr_m^{\am+1}\right)^{\frac{\am}{\am+1}} |J\times K|^{\frac{1}{\am +1}}, 
\end{equation} 
since  $\frac{1}{\am +1}\to 0$ ($\alpha_{m}\to +\infty$),  then using \eqref{vr-m+1}
\begin{equation}
\label{d1}
\int_{J\times K} \vr_m^\am \dx\dt \leq C({\rm data}, n) |J\times K|.
\end{equation}
Thus we can conclude that 
\begin{equation}
\pi_m(\vr_m)\rightharpoonup \pi \ \mbox{ as } m\to +\infty \mbox{ in }  L^1(J\times K),
\end{equation}
 for any compact set $J\times K\subset (0,T)\times \Omega$ at least for a chosen subsequence.

\subsubsection{The limit as $m\to +\infty$ in the equation of continuity and momentum equation}
We are in position to take the limit as $m\to +\infty$. For any $n\in \mathbb{N}$ it holds
\begin{equation} \label{MomEq-n}
\begin{split}
\left[ \intO{ \vr \vu \cdot \bfphi } \right]_{t=0}^{t = \tau} &= 
\int_0^\tau \intO{ \Big[ \vr \vu \cdot \partial_t \bfphi + \vr \vu \otimes \vu : \Grad \bfphi 
+ {\pi} \Div \bfphi - \mathbb{S}  : \Grad \bfphi \Big] }
\end{split}
\end{equation}
for any $\bfphi \in C^1([0,T]; X_n)$.
 Passing into to limit in \eqref{E16} we obtain the following 
\begin{equation} \label{continuity-n}
\begin{split}
\left[ \intO{ \vr \varphi } \right]_{t = 0}^{t = \tau}&= 
\int_0^\tau \intO{ \Big[ \vr \partial_t \varphi + \vr \vu \cdot \Grad \varphi  \Big] }
\dt \\ 
&- \int_0^\tau \int_{\Gamma_{\rm out}} \varphi \vr \vu_B \cdot \vc{n} \ {\rm d} S_x \dt - 
\int_0^\tau \int_{\Gamma_{\rm in}} \varphi \vr_B \vu_B \cdot \vc{n}  \ {\rm d}S_x \dt,
\end{split}
\end{equation}
for any $\varphi \in C^1([0,T] \times \Ov{\Omega})$ and for any $n\in \mathbb{N}$.

\subsubsection{The limit as $m\to +\infty$ in the approximate energy balance}
Since  $\vr_m^\am$ is uniformly bounded  
and $\am \to +\infty$ as $m\to +\infty$  it results that  
\begin{equation}
\Pi_m(\vr_m)= \frac{1}{\am-1} \vr_m^\am \to 0 \mbox{  in } L^1((0,T)\times\Omega).
\end{equation}
Then employing the convergences \eqref{u-m-infty} -- \eqref{rho-u-m} and finally the weak lower semicontinuity of the potentials $F$ and $F^*$ we may perform the limit as $m\to +\infty$ in the energy balance \eqref{En-ep} obtaining that for any $n\in \mathbb{N}$ it holds
\begin{equation} \label{Energy-Galerkin}
\begin{split}
\left[ \intO{\frac{1}{2} \vr |\vu - \vu_B|^2 } \right]_{t = 0}^{ t = \tau} &+ 
\int_0^\tau \intO{\left[ F(\Ds \vu) + F^*(\mathbb{S}) \right] } \dt 
 +\int_0^\tau \intO{ \pi \Div \vu_B}\dt  \\
&\leq- 
\int_0^\tau \intO{ \left[ \vr \vu \otimes \vu \right] :  \Grad \vu_B } \dt \\ &+ \int_0^\tau \intO{ {\vr} \vu  \cdot \vu_B \cdot \Grad \vu_B  } 
\dt
+ \int_0^\tau \intO{ \mathbb{S} : \Grad \vu_B } \dt. 
\end{split}
\end{equation}

\subsubsection{The free boundary condition}
In order to conclude this step we need to recover the free boundary condition. From \eqref{positivepart} as $0\leq \vr\leq 1$ we already know that
 $\vr\pi\leq \pi,$ (for details see \cite{LiMa}, 
 then in order to prove that $\vr\pi=\pi$ we need to show that $\vr\pi\geq \pi$. To this aim we shall use  in the approximating momentum equation \eqref{MM} the test function 
$$\bfphi_m(t,x)= \eta(t)\, \phi(x) \Grad \Del^{-1}(\phi(x) \vr_m(t,x)), \ \eta\in C^\infty_{\rm c}(0,T), \  \phi\in C^\infty_{\rm c}(\Omega) $$
and in the limiting momentum equation \eqref{MomEq-n} the test function 
$$\bfphi(t,x)= \eta(t)\, \phi(x) \Grad \Del^{-1}(\phi(x) \vr(t,x)), \ \eta\in C^\infty_{\rm c}(0,T), \  \phi\in C^\infty_{\rm c}(\Omega). $$
We subtract both identities and by performing the limit as $m\to +\infty$, taking into account the convergences \eqref{u-m-infty}--\eqref{rho-u-m}  and that since $n$ is fixed in particular it holds that $\vr_m\Div\vu_m \to\vr\Div\vu$ weakly-* in $L^\infty((0,T)\times\Omega)$
we obtain
\begin{equation}
\int_0^T\eta \intO{ \phi^2\,\pi\vr}-\lim_{m\to +\infty}\int_0^T\eta \intO{ \phi^2\,\pi_m(\vr_m)\vr_m} =0.
\end{equation}
Therefore 
\begin{equation}
\Ov{\vr_m^{\am +1}} =\vr\pi,
\end{equation}
where $\Ov{\vr_m^{\am +1}}$ is the weak limit of $\vr_{m}^{\alpha_{m}+1}$.
Since $\vr_m\geq 0$ it holds that for any $\varepsilon>0$ there exists $k$ such that for any $m\geq k$ 
$$ \vr_m^{\am +1} \geq \vr_m^\am -\ep,$$
passing to the weak limit we deduce that 
$$\vr \pi=\Ov{\vr_m^{\am +1}} \geq \pi-\ep,$$
then taking the limit as $\ep\to 0$ we get $\vr\pi\geq \pi$ and finally
\begin{equation}
\label{id=press}
\vr\pi=\pi.
\end{equation}

\subsubsection{The divergence free condition}
In order to conclude the limit passage as $m\to +\infty$, 
we prove that the divergence free condition $\Div \vu=0$ is satisfied a.e.  in $\{ \vr=1\}$. We show that it follows from a compatibility condition between equation \eqref{c-e} and conditions \eqref{bounded}, \eqref{r=1}.

\bLemma{Eq}{Let $\vu\in L^q(0,T; W^{1,q}(\Omega; \mathbb{R}^3))$, such that $\vu-\vu_B\in  L^q(0,T; W_0^{1,q}(\Omega; \mathbb{R}^3))$ and $\vr \in L^2((0,T)\times\Omega)$ such that
$$  \partial_t \vr + \Div (\vr \vu) = 0  \ \mbox{ in } (0, T)\times\Omega, \  \vr(0, \cdot) = \vr_0,$$
$$\vr|_{\Gamma_{in}} = \vr_B,\ 0<\vr_B\leq 1, \ 
\Gamma_{in} = \left\{ x \in \partial \Omega \ \Big| \ \vu_B \cdot \vc{n} < 0 \right\},$$
then the following two assertions are equivalent
\begin{itemize}
\item[(i)] $\Div \vu=0$ a.e. on $ \{\vr\geq 1\}$ and $ 0< \vr_0\leq 1$,
\item[(ii)] $ 0< \vr\leq 1.$
\end{itemize}
}
\eL
\noindent
{\bf Proof. } The proof follows the arguments in \cite[Lemma 2.1]{LiMa} and \cite[Theorem 11.36]{FeNo-book}, with small changes when proving $(i) \Rightarrow (ii)$ due to the boundary condition $\vu|_{\partial\Omega}=\vu_B$. For the sake of clarity we perform here the proof of such implication  $(i) \Rightarrow (ii)$. 
From the regularization lemma stated by P.L.~Lions \cite{PL-L} it holds
$$\partial_t\beta(\vr) +\Div(\beta(\vr)\vu) = (\beta(\vr) -\vr\beta'(\vr))\Div \vu $$
for any $C^1$ function $\beta:\mathbb{R}\to\mathbb{R}$ such that $|\beta(t)|\leq C(1 +t)$ with $C$ positive constant.  
 Let $\beta$ be the following function
\begin{equation*}
\beta(t)=\begin{cases}
0 \  &\mbox{ if } t\leq0,\\
t \ &\mbox{ if } 0\leq t\leq 1,\\
 1 \ &\mbox{ if } t\geq 1,
\end{cases}
\end{equation*}
then employing assumption $(i)$ we get 
\begin{equation}
\partial_t\beta(\vr) +\Div(\beta(\vr)\vu) = 0
\end{equation}
To formalize this process it is enough to  consider a sequence of $C^1$ functions $\beta_n$ converging to $\beta$ pointwise and in $L^2$. Setting
 $d=\beta(\vr)-\vr$, we have that $d$ solves the equation
 $\partial_t d + \Div( d\vu)=0$.  Employing once more the P.L.~Lions regularization  \cite{PL-L} we get that $|d|$ solves the same equation. Integrating in space we get 
 that 
 $$\intO{\partial_t(|d|)}= -\int_{\Omega\setminus \Gamma_{in}} |d(\vr)| \vu_B\cdot\vc{n} {\rm d}\sigma_x\leq 0.$$ 
Integrating in time and observing that $|d||_{t=0}=0$ we deduce
 $$\intO{|d|(t)}=0,$$ 
 hence $\beta(\vr)\leq \vr$ which yields $0\leq \vr\leq 1$,

\subsection{Limit in the Galerkin approximations}
Here we perform the last level approximation limit $n\to +\infty$. As a consequence  the uniform bounds on the velocity $\vu_{n}$ does not hold anymore since we are not  in a finite dimensional space. This fact entails that we have to perform some more refined estimates and we loose the strong convergence of the momentum. As we will see later on this fact gives rise to the dissipative measures.

\subsubsection{Uniform estimates} 
In order to recover the uniform bounds from the energy inequality \eqref{Energy-Galerkin} we need to employ the  Gronwall inequality. We see that at this level it is necessary to have information on the sign 
of  the term 
$$ \int_0^\tau \intO{ \pi \Div \vu_B}\dt,$$
where with the same notation $\vu_B$ we denoted an extension inside $\Omega$ of the boundary velocity. 
We recall the assumption \eqref{bdu-bis} on the boundary data
\begin{equation}\label{assumption-ub}
\vu_B\in C^2(\partial\Omega;\mathbb{R}^d),  \int_{\partial\Omega} \vu_B\cdot\vn \geq 0, \quad  \Omega \mbox{ is a } C^2 \mbox{ bounded domain,}
\end{equation}
and  following Choe-Novotn\'{y}-Yang \cite[Lemma 5.2]{ChNoYa} and Pokorny-Wr\'{o}blewska-Zatorska \cite[Lemma 2.5]{PoWrZa} we claim that  there exists an extension (that we still denote $\vu_B$) such that 
\begin{equation}
\Div \vu_B \geq 0 \mbox{ a.e. in } \Omega.
\end{equation}
Now, thanks also to the fact that $\pi_n\geq 0$, we are in position to claim that 
\begin{equation}
 \int_0^\tau \intO{ \pi_n \Div \vu_B}\dt \geq 0.
\end{equation}
At this point from \eqref{Energy-Galerkin}  employing the regularity of the boundary velocity $\vu_B$ and the growth assumption \eqref{S5a} on $F$ ,  the application of the Gronwall inequality provides the following uniform estimates as $n\to +\infty$ 
\begin{align}\label{vu-n}
&\sup_{n\in \mathbb{N}}\ \sup_{t\in [0,T]}\intO{\frac{1}{2} \vr_n |\vu_n - \vu_B|^2 } <+\infty,\\
&\sup_{n\in \mathbb{N}} \|\vu_n\|_{L^q(0, T; W^{1,q}(\Omega; \mathbb{R}^d))}<+\infty,\\
&\sup_{n\in \mathbb{N}} \intTO{F^*(\mathbb{S}_n)}<+\infty.\label{FstarS}
\end{align}
Moreover, from \eqref{rho-limitata} and \eqref{rho-bordo-limitata} we know that 
\begin{align}
&0\leq \vr_n\leq 1 \mbox{ a.e. in } (0,T)\times\Omega,\\
&0\leq \vr_n\leq 1 \mbox{ a.e. in } (0,T)\times\Gamma_{\rm out} \mbox{ with respect to the measure on the boundary } |\vu_B\cdot\vn|  \D S_x.
\end{align}
In particular the sequence $\vr_n$ is uniformly bounded in $L^p((0,T)\times\Omega)\cap L^p((0,T)\times \partial\Omega)$ for any $1\leq p<+\infty$. Therefore, up to a subsequence and recalling the de la Vall\`{e}-Poussin criterion, we have 
\begin{align}
&\vu_n \to \vu \mbox{ weakly in } L^q(0, T; W^{1,q}(\Omega; \mathbb{R}^d),\label{u-q}\\
&\mathbb{S}_n \to \mathbb{S} \mbox{ weakly in } L^1((0,T)\times\Omega; \mathbb{R}^{d\times d}), \label{S-1}\\
&\vr_n \to \vr \mbox{ in } C_{{\rm weak}}([0, T]; L^p(\Omega)), \quad 1\leq p<+\infty,\label{rho-p-omega}\\
&\vr_n|_{\Gamma_{\rm out}}\to\vr  \mbox{ weakly-*  in } L^\infty((0,T); L^p(\Gamma_{\rm out}; |\vu_B\cdot\vn|  \D S_x)),  \quad 1\leq p<+\infty,\label{rho-p-boundary}\\
&\vr_n\vu_n \to \vm  \mbox{ weakly-*  in } L^\infty(0,T; L^{\frac{2p}{p+1}}(\Omega; \mathbb{R}^d)) \quad 1\leq p<+\infty.\label{m}
\end{align}

In order to prove the convergence of the sequence of measures $\pi_{n}$, 
let us choose in \eqref{MomEq-n} the following test function 
$$ \bfphi=\eta(t)\mathcal{B}[\xi]$$
where $\mathcal{B}$ is the Bogovskii operator, $\xi$ is a $C_c^\infty(\Omega)$ function such that $\intO{\xi}=0$, and  $\eta$ is in $C^1_{\rm c}(0,T)$. Then we have 
\begin{equation}
\int_0^T\intO{\eta \pi_n \xi}\dt=\int_0^T \intO{ \Big[ \eta' \,\vr_n \vu_n \cdot \mathcal{B}[\xi] + \eta\left[ \vr_n \vu_n \otimes \vu_n -  \mathbb{S}_n \right]: \Grad \mathcal{B}[\xi]\Big] }\dt.
\end{equation}
Employing the uniform estimates \eqref{vu-n}--\eqref{FstarS} we get the uniform boundedness of the right hand-side, thus up to a subsequence we  deduce that
\begin{equation}
\pi_n\to \pi \mbox{ weakly-* in } \mathcal{M}^+(I\times K) \mbox{ for any compact } I\times K\subset (0,T)\times \Omega.
\end{equation}

\subsubsection{The free boundary condition}
For any $n\in\mathbb{N}$, $\eta\in C_c^\infty(0,T)$, $\varphi\in C_c^\infty(\Omega)$ from \eqref{id=press} it holds 
$$\int_0^T \!\!\intO{\eta \varphi \vr_n\pi_n}\dt= \int_0^T\!\! \intO{\eta \varphi \pi_n}\dt$$
taking the limit as $n\to+\infty$ and employing the convergences showed in the previous subsection we get 
$$\int_0^T \!\!\intO{\eta \varphi \vr\pi}\dt= \int_0^T\!\!\intO{\eta \varphi \pi}.\dt$$
Therefore 
\begin{equation}
(\vr-1)\pi =0 \mbox{ a.e. in } (0,T)\times\Omega.
\end{equation}

\subsubsection{Limit in the continuity equation}
By virtue of the convergences \eqref{rho-p-omega}--\eqref{m} we are in position to take the limit as $n\to +\infty$ in \eqref{continuity-n} and obtain the following 
\begin{equation} \label{continuity-finale}
\begin{split}
\left[ \intO{ \vr \varphi } \right]_{t = 0}^{t = T}&= 
\int_0^T \intO{ \Big[ \vr \partial_t \varphi + \vr \vu \cdot \Grad \varphi  \Big] }
\dt \\ 
&- \int_0^T \int_{\Gamma_{\rm out}} \varphi \vr \vu_B \cdot \vc{n} \ {\rm d} S_x \dt - 
\int_0^T \int_{\Gamma_{\rm in}} \varphi \vr_B \vu_B \cdot \vc{n}  \ {\rm d}S_x \dt,
\end{split}
\end{equation}
for any $\varphi \in C^1([0,T] \times \Ov{\Omega})$. Also $0\leq \vr\leq 1$ and $\vr(0,\cdot)=\vr_0$.

\subsubsection{Limit in the momentum equation and in the energy inequality}
In order to pass into the limit in the momentum equation \eqref{MomEq-n} we observe that from \eqref{vu-n} we have 
\begin{equation}
\begin{split}
\vr_n \vu_n \otimes \vu_n&=
1_{\vr_n > 0} \frac{\vm_n \otimes \vm_n}{\vr_n} 
\to \Ov{ 1_{\vr > 0} \frac{\vm \otimes \vm}{\vr} } \mbox{weakly-(*) in}\ L^\infty(0,T; \mathcal{M}(\Ov{\Omega}; R^{d \times d}_{\rm sym})),
\end{split}
\end{equation}
and
\begin{equation}\label{kinetic}
\begin{split}
\frac{1}{2} \vr_n |\vu_n|^2 =
\frac{1}{2} \frac{|\vm_n|^2}{\vr_n}  \to \Ov{  \frac{1}{2} \frac{|\vm|^2}{\vr} }
\mbox{weakly-(*) in}\ L^\infty(0,T; \mathcal{M}(\Ov{\Omega})).
\end{split}
\end{equation}
We set
\[
\begin{split}
\mathfrak{R} = \Ov{ 1_{\vr > 0} \frac{\vm \otimes \vm}{\vr} } -
1_{\vr > 0} \frac{\vm \otimes \vm}{\vr} = \Ov{ 1_{\vr > 0} \frac{\vm \otimes \vm}{\vr} } -  \vr \vu \otimes \vu,
\end{split}
\]
and
\[
\mathfrak{E} =
\Ov{ \frac{1}{2} \frac{|\vm|^2}{\vr}  } -
\frac{1}{2} \frac{|\vm|^2}{\vr}=  \frac{1}{2} \frac{|\vm|^2}{\vr}  
-  \frac{1}{2} \vr |\vu|^2 
\]
noting the relation
\begin{equation}\label{E-R}
\underline{d} \mathfrak{E} \leq {\rm tr}[ \mathfrak{R}] \leq \Ov{d} \mathfrak{E},
\ \mbox{where}\ 0 < \underline{d} \leq \Ov{d}, \ \underline{d} = \underline{d} (d).
\end{equation}
We observe that
\[
\mathfrak{R} \in L^\infty(0,T; \mathcal{M}^+ (\Ov{\Omega}; R^{d \times d}_{\rm sym})).
\]
Indeed we compute
\[
\left[ \Ov{ 1_{\vr > 0} \frac{ \vm \otimes \vm }{\vr} } - 1_{\vr > 0} \frac{\vm \otimes \vm}{\vr} \right]
:(\xi \otimes \xi) = \left[ \Ov{ \frac{ |\vm \cdot \xi|^2 }{\vr} } - \frac{ |\vm \cdot \xi|^2 }{\vr} \right] \geq 0
\ \mbox{for any}\ \xi \in R^d,
\]
{where the most right inequality follows from the convexity of the lower semicontinuous function
$$
[\vr, \vm] \mapsto \left\{\begin{array}{c}
 \frac{|\vm \cdot \xi|^2}{\vr}\;\mbox{if $\vr>0$},
 \\
 0\;\mbox{if $\vr=0$, $\vm=0$},\\
 \infty\;\mbox{otherwise}.
                   \end{array}\right.
$$
}

Now, taking the limit in \eqref{MomEq-n} we get 
\begin{equation} \label{Momentum-finale}
\begin{split}
\left[ \intO{ \vr \vu \cdot \bfphi } \right]_{t=0}^{t = T} &= 
\int_0^T \intO{ \Big[ \vr \vu \cdot \partial_t \bfphi + \vr \vu \otimes \vu : \Grad \bfphi 
 - \mathbb{S}  : \Grad \bfphi \Big] } \\
 &+\langle \pi, \Div\bfphi\rangle_{(\mathcal{M}((0,T)\times\Omega), C((0,T)\times\Omega))} + \int_0^T\int_\Omega\Grad\bfphi:\D\mathfrak{R}(t)\dt
\end{split}
\end{equation}
for any $\bfphi \in C_c^1((0,T)\times\Omega; \mathbb{R}^d)$.

Then employing the convergences \eqref{u-q}, \eqref{S-1}, \eqref{kinetic} and relation \eqref{E-R}, finally the weak lower semicontinuity of the potentials $F$ and $F^*$ we may perform the limit as $n\to +\infty$ in the energy balance \eqref{Energy-Galerkin} obtaining that it holds
\begin{equation} \label{Energy-final}
\begin{split}
\left[ \intO{\frac{1}{2} \vr |\vu - \vu_B|^2 } \right]_{t = 0}^{ t = T} &+ 
\int_0^T \intO{\left[ F(\Ds \vu) + F^*(\mathbb{S}) \right] } \dt 
+\langle \pi, \Div\vu_B\rangle_{(\mathcal{M}((0,T)\times\Omega), C((0,T)\times\Omega))}\\
   + \frac{1}{\Ov{d}} \int_{\Ov{\Omega}}  \D \ {\rm tr}[\mathfrak{R}] (T)
&\leq -\int_0^\tau \int_{\Ov{\Omega}} \Grad \vu_B : \D \ \mathfrak{R}(t) \dt - 
\int_0^T \intO{ \left[ \vr \vu \otimes \vu \right] :  \Grad \vu_B } \dt \\ &+ \int_0^\tau \intO{ {\vr} \vu  \cdot \vu_B \cdot \Grad \vu_B  } 
\dt
+ \int_0^T \intO{ \mathbb{S} : \Grad \vu_B } \dt
\end{split}
\end{equation}
with $\Ov{d}$ a constant depending on the dimension $d$. This final step concludes the proof of the Theorem \ref{T1}.

\begin{Remark}
We conclude this paper by describing two related models for the problem \eqref{problem}, where the condition \eqref{pi=0} can be generalized.

\underline{\emph{General fluid pressure law}}
The free-boundary conditions \eqref{pi=0} can be  generalised by including a general fluid pressure (baratropic pressure), namely 

$$\pi=p(\vr)\  \mbox{ on } \{\vr<1\}, \qquad  \pi\geq p(1) \mbox{  on } \{\vr=1\}$$
or equivalently
$$ \vr (\pi - p(\vr)) = \pi - p(\vr).$$

The fluid behaviour in this case is that of a barotropic fluid in the region $ \{ \vr < 1\}$. This generalization requires only some  changes in the energy estimates which can be treated with similar techniques.

\underline{\emph{Congestion constraints.}}
With a similar  analysis  we  can  also deal  with a non-homogeneous congestions constraints, i.e. a non homogeneous threshold for the pressure. In this case \eqref{pi=0} has the form
$$\pi=0\  \mbox{ on } \{\vr^{*}<1\}, \qquad  \pi\geq 0 \mbox{  on } \{\vr=\vr^{*}\}$$
In this case we need to change the approximating pressure. Instead of \eqref{pi-m} we take an approximating pressure of the form
$$\pi_m = \left(\frac{\vr}{\vr^{*}}\right)^{\alpha_m}.$$
\end{Remark}

\section{Compatibility with strong solutions.}
\label{compatibility}
We conclude this paper by showing the compatibility between dissipative solutions enjoying additional regularity and classical solutions. 
Let us assume that a dissipative solution (whose existence is proved in the previous sections) enjoys the following regularity
\begin{equation}
\vu\in C^1([0,T]\times \Ov{\Omega}; \mathbb{R}^d), \ \vr\in C^1([0,T]\times \Ov{\Omega}), \inf_{(0,T)\times\Omega}\vr >0, \ \pi\in C([0,T]\times\Ov{\Omega}),
\end{equation}
then $[\vr,\vu,\pi]$ is a classical solution, meaning that the measure $\mathfrak{R}=0$.
 Let us consider $\vu-\vu_B$ as test function in the momentum equation \eqref{Momentum-finale}
\begin{equation}\begin{split}
&\left[ \intO{ \vr (\vu-\vu_B) \cdot (\vu-\vu_B) } \right]_{t=0}^{t = \tau} + \left[ \intO{ \vr \vu_B \cdot (\vu-\vu_B) } \right]_{t=0}^{t = \tau}= \\
&\int_0^\tau \intO{ \Big[\frac{1}{2} \vr \partial_t  |\vu-\vu_B|^2 + \vr \vu \otimes \vu : \Grad  (\vu-\vu_B) 
 - \mathbb{S}  : \Grad  (\vu-\vu_B)  \Big] } \\
 &+\int_0^\tau \intO{ \vr \vu_B \cdot \partial_t  (\vu-\vu_B)} +\langle \pi, \Div (\vu-\vu_B) \rangle_{(\mathcal{M}((0,T)\times\Omega), C((0,T)\times\Omega))} \\
 &+ \int_0^\tau\int_\Omega\Grad (\vu-\vu_B) :\D\mathfrak{R}(t)\dt.
\end{split}\end{equation}
By using the continuity equation \eqref{continuity-finale}, we get
\begin{equation}\label{comp}\begin{split}
&\frac{1}{2}\left[ \intO{ \vr |\vu-\vu_B|^2  } \right]_{t=0}^{t = \tau} =- \frac{1}{2} \int_0^\tau \intO{ \vr \vu\cdot \Grad  |\vu-\vu_B|^2}\dt\\
&+\int_0^\tau \intO{\Big[ \vr \vu \otimes \vu : \Grad  (\vu-\vu_B)  - \mathbb{S}  : \Grad  (\vu-\vu_B)  \Big] }\dt \\
&- \int_0^\tau \intO{ \vr\vu\cdot \vu_B \cdot \Grad  (\vu-\vu_B)}\dt  - \int_0^\tau \intO{ \vr\vu\cdot \Grad\vu_B \cdot (\vu-\vu_B)}\dt\\
 &+\langle \pi, \Div (\vu-\vu_B) \rangle_{(\mathcal{M}((0,T)\times\Omega), C((0,T)\times\Omega))}+ \int_0^\tau\int_\Omega\Grad (\vu-\vu_B) :\D\mathfrak{R}(t)\dt.
\end{split}\end{equation}
Subtracting  \eqref{comp} from \eqref{Energy-final} and reorganizing the various terms  it yields that
\begin{equation}\label{comp2}\begin{split}
 \frac{1}{\Ov{d}} \int_{\Ov{\Omega}}  \D \ {\rm tr}[\mathfrak{R}] (\tau)+ &\int_0^\tau \intO{\left[ F(\Ds \vu) + F^*(\mathbb{S}) - \mathbb{S}  : \Grad  \vu \right] } \dt \\ &+\langle \pi, \Div \vu\rangle_{(\mathcal{M}((0,T)\times\Omega), C((0,T)\times\Omega))}
 \leq  \int_0^\tau\int_\Omega\Grad \vu :\D\mathfrak{R}(t)\dt.
\end{split}\end{equation}
Observing that from free boundary conditions \eqref{r=1} and \eqref{pi=0} it follows that
$$\langle \pi, \Div \vu\rangle_{(\mathcal{M}((0,T)\times\Omega), C((0,T)\times\Omega))}=0,$$
and from the Fenchel-Young inequality it holds
$$F(\Ds \vu) + F^*(\mathbb{S}) - \mathbb{S}  : \Grad  \vu\geq 0,$$
we can employ the Gronwall lemma in \eqref{comp2} yielding that $\mathfrak{R}=0$, then going back into \eqref{comp2} we deduce that $\mathbb{S}\in \partial F(\mathbb{D}\vu)$.

\subsection*{Acknowledgments}

The authors gratefully acknowledge the partial support by the Gruppo
Na\-zio\-na\-le per l’Analisi Matematica, la Probabilit\`a e le loro
Applicazioni (GNAMPA) of the Istituto Nazionale di Alta Matematica
(INdAM), and by the PRIN 2020- Project N. 20204NT8W4 ``Nonlinear evolution PDEs, fluid
dynamics and transport equations: theoretical foundations and
applications''. 
A.A. is member of the Italian National Group for the Mathematical Physics (GNFM) of the Italian National Institute for the High Mathematics (INdAM) for the year 2024. 
\subsection*{Declarations}

\bigskip
\noindent
\textbf{Data Availability.} The authors declare that data sharing is not applicable to this article as no datasets were generated or analysed.

\bigskip
\noindent
\textbf{Conflict of interest.} The authors declare that they have no conflict of interest.

\bibliographystyle{amsplain}

\providecommand{\bysame}{\leavevmode\hbox to3em{\hrulefill}\thinspace}
\providecommand{\MR}{\relax\ifhmode\unskip\space\fi MR }
\providecommand{\MRhref}[2]{%
  \href{http://www.ams.org/mathscinet-getitem?mr=#1}{#2}
}
\providecommand{\href}[2]{#2}

\end{document}